\documentclass [reqno,10pt,oneside]{amsart}

\usepackage{amsthm}
\usepackage{amscd}
\usepackage{amsmath,amssymb}
\usepackage{algorithm}
\usepackage[noend]{algorithmic}

\theoremstyle {plain}
\newtheorem {thm}{Theorem}[section]
\newtheorem {prop}[thm]{Proposition}

\theoremstyle {definition}

\theoremstyle {remark}
\newtheorem {rem}[thm]{Remark}

\newtheorem {exmp}[thm]{Example}

\newcommand{\A}{{\mathcal A}}
\newcommand{\E}{{\mathcal E}}

\newcommand{\V}{{\mathcal V}}
\newcommand{\W}{{\mathcal W}}
\newcommand{\C}{{\mathbb C}}

\newcommand{\N}{{\mathbb N}}
\newcommand{\Q}{{\mathbb Q}}
\newcommand{\R}{{\mathbb R}}
\newcommand{\Z}{{\mathbb Z}}

\newcommand{\singular}{{\sc Singular }}

\begin{document}

\bibliographystyle{alpha}

\title{About the Computation of the Signature of Surface Singularities $z^N+g(x,y)=0$}

\author{Muhammad Ahsan Banyamin}
\address{Muhammad Ahsan Banyamin\\ Abdus Salam School of Mathematical Sciences\\ GC University\\ 
Lahore\\ 68-B\\ New Muslim Town\\ Lahore 54600\\ Pakistan}
\email{ahsanbanyamin@gmail.com}

\author{Gerhard Pfister}
\address{Gerhard Pfister\\ Department of Mathematics\\ University of Kaiserslautern\\
Erwin-Schr\"odinger-Str.\\ 67663 Kaiserslautern\\ Germany}
\email{pfister@mathematik.uni-kl.de}
\urladdr{http://www.mathematik.uni-kl.de/$\sim$pfister} 

\author{Stefan Steidel}
\address{Stefan Steidel\\ Department of Mathematics\\ University of Kaiserslautern\\ 
Erwin-Schr\"odinger-Str.\\ 67663 Kaiserslautern\\ Germany}
\email{steidel@mathematik.uni-kl.de}
\urladdr{http://www.mathematik.uni-kl.de/$\sim$steidel} 


\thanks{Part of the work was done at ASSMS, GCU Lahore -- Pakistan.}

\date{\today}

\maketitle

\begin{abstract}
In this article we describe our experiences with a parallel \textsc{Singular} -- implementation
of the signature of a surface singularity defined by $z^N+g(x,y)=0$.
\end{abstract}

\section{Introduction}

Let $g \in \C[x,y]$ define an isolated curve singularity at $0 \in \C^2$ and $f:=z^N+g(x,y)$ for $N \geq 2$.
The zero--set $V:=V(f) \subseteq \C^3$ of $f$ has an isolated singularity at $0$.
For a small $\varepsilon > 0$ let $V_\varepsilon := V(f-\varepsilon) \subseteq \C^3$ be the Milnor
fibre of $(V,0)$ and $s: H_2(V_\varepsilon,\R) \times H_2(V_\varepsilon,\R) \longrightarrow \R$
be the intersection form (cf. \cite{AGZV}, \cite{M}, \cite{N95CM}, \cite{St1}). $H_2(V_\varepsilon,\R)$ is a 
$\mu$--dimensional $\R$--vector space, $\mu$ the Milnor number of $(V,0)$ (cf. \cite{AGZV},
\cite{JP}, \cite{KN}, \cite{K}), and $s$ is a symmetric bilinear form. 
Let $\sigma(f)$ be the signature of $s$, called the signature of the surface singularity $(V,0)$.
Formulas to compute the signature are given by N\'emethi (cf. \cite{N98}, \cite{N99}) and 
van Doorn, Steenbrink (cf. \cite{DS}).
We have implemented three approaches in \textsc{Singular} (cf. \cite{DGPS}, \cite{GP}) using Puiseux 
expansions, the resolution of singularities and the spectral pairs of the singularity.
In section \ref{secTheory} we will explain the different formulas to compute the signature, 
and finally we give examples and timings of our implementation in section \ref{secExTime}.

\section{The signature of $(V,0)$ in terms of $N$ and $g$} \label{secTheory}

The proofs of the following formulas (cf. Propositions \ref{propPuiseux}, \ref{propBrieskorn},
\ref{propResolution} and \ref{propSpectralPairs}) can be found in the corresponding papers 
of N\'emethi (cf. \cite{N98}, \cite{N99}). 

\subsection{Approach 1: Puiseux Pairs} \label{subsecPuiseuxPairs}

For the first approach assume that $(V(g),0) \subseteq (\C^2,0)$ is the germ of an
irreducible curve singularity. Let $(m_1,n_1), \ldots, (m_\ell,n_\ell)$ be the \emph{Puiseux pairs}
of $g$ and define a sequence $\{a_i\}_{i=1,\ldots,\ell}$ by $$a_1 = m_1 \quad \text{ and }
\quad a_{i+1} = m_{i+1} - n_{i+1} \cdot (m_i - n_i \cdot a_i).$$
Moreover, we set $d_\ell = 1$ and $d_i = \gcd(N,n_{i+1} \cdots n_\ell)$ for $1 \leq i < \ell$.

\begin{prop} \label{propPuiseux}
$\sigma(f) = \sum_{i=1}^\ell d_i \cdot \sigma\left(x^{a_i} + y^{n_i} + z^{N/d_i}\right)$.
\end{prop}

The signature of a Brieskorn polynomial $x^{c_1} + y^{c_2} + z^{c_3}$ can be computed
combinatorially. Let $S_t = \#\{(k_1,k_2,k_3) \in \Z^3 \mid 1 \leq k_j \leq a_j -1, 
t < \sum_{j=1}^3 \frac{k_j}{a_j} < t+1\}$ for $t \in \N_0$.

\begin{prop} \label{propBrieskorn}
$\sigma(x^{c_1} + y^{c_2} + z^{c_3}) = S_0-S_1+S_2$.
\end{prop}

\begin{rem}
The \textsc{Singular} implementation of the first approach bases on the procedure \texttt{invariants}
(cf. library \texttt{hnoether.lib}) to obtain the Puiseux pairs of $g$. The command 
\texttt{list L = invariants(g);}\footnote{The command \texttt{list L = invariants(g);} returns a list \texttt{L} 
of the following format: \texttt{L[1]}: characteristic exponents, \texttt{L[2]}: generators of the semigroup, 
\texttt{L[3]}: first components of Puiseux pairs, \texttt{L[4]}: second components of Puiseux pairs, \texttt{L[5]}: 
degree of the conductor, \texttt{L[6]}: sequence of multiplicities.} provides a list \texttt{L} in which the third 
and fourth entry contain the necessary data for our application. The required combinatorics of Proposition 
\ref{propPuiseux} resp. \ref{propBrieskorn} had to be implemented in spring 2011.
\end{rem}

\subsection{Approach 2: Resolution} \label{subsecResolution}

In the following  $(V(g),0)$ doesn't need to be irreducible. For the second approach we use the 
\emph{resolution} of the singularity $(V(g),0)$. Let $\V = \W \coprod \A$ be the vertices of the
resolution graph, $\W$ the vertices corresponding to the exceptional divisors and $\A$ the vertices 
corresponding to the resolved branches. Let $\E = \{ (v_1,v_2) \mid v_1,v_2 \in \V\}$ be the set 
of edges of the resolution graph.
Let $\{m_v\}_{v \in \W}$ be the sequence of total multiplicities and set $m_a = 1$ if $a \in \A$. 
For $w \in \W$ let $M_w = \gcd(m_v \mid v \in \V_w \cup \{w\})$ and for $e = (v_1,v_2) \in \E$ let 
$m_e = \gcd(m_{v_1}, m_{v_2})$.

\begin{prop} \label{propResolution}
$\sigma(z^N+g) = \eta(g,N) - N \cdot \eta(g,1)$ and \footnote{For $x \in \Q$ we denote by $\{x\}$
the fractional part of $x$ and $$((x)) = \left\{ \begin{array}{rl} \{x\}-\frac 12, & \text{if } x \notin \Z \\
0, & \text{if } x \in \Z \end{array} \right..$$ The definition of the \emph{eta--invariant} $\eta(f,K)$
can be found in \cite{N95T}.}
\begin{align*}
\eta(g,K) = & \#(\A) - 1 + \sum_{e \in \E} \big( \gcd(K,m_e)-1\big) - \sum_{w \in \W} \big( \gcd(K,M_w)
                        -1 \big) \\
                    & + 4 \cdot \sum_{\substack{w \in \W \\ \text{$w$ \emph{has more}} \\ \text{\emph{than $2$ 
                        edges}}}} \sum_{\substack{v \in \V \\ (v,w) \in \E}} \sum_{k=1}^{m_w} ((\frac{k \cdot m_v}
                        {m_w})) \cdot ((\frac{k \cdot K}{m_w})).
\end{align*}
Moreover, it holds $\sum_{e \in \E} \big( \gcd(K,m_e)-1\big) = \sum_{w \in \W} \big( \gcd(K,M_w)
-1 \big)$ if $(V(g),0)$ is irreducible.
\end{prop}

\begin{rem}
The \textsc{Singular} implementation of the second approach bases on the procedure 
\texttt{totalmultiplicities} (cf. library \texttt{alexpoly.lib}) to obtain necessary information about the 
resolution of $(V(g),0)$. The command \texttt{list L = totalmultiplicities(g);} provides a list \texttt{L}
of the following format: \texttt{L[1]}: incidence matrix of the resolution graph, \texttt{L[2]}: list of the 
sequences of the total multiplicities corresponding to the branches, \texttt{L[3]}: list of the multiplicity 
sequences of the branches. The required combinatorics of Proposition \ref{propResolution} had to be 
implemented in spring 2011.
\end{rem}

\subsection{Approach 3: Spectral Pairs} \label{subsecSpectralPairs}

The third approach uses the \emph{spectral pairs}\footnote{A definition of \emph{spectral pairs}
can be found in \cite{K}.} of the singularity $(V(g),0)$.
Therefore let $$Spp(g) = \sum_{(\alpha,w)} h^{1+[-\alpha], w + s_\alpha -1 - [-\alpha]}_{\exp(-2\pi i\alpha)}
\cdot (\alpha,w)$$ represent the spectral pairs where $s_\alpha = 0$ if $\alpha \notin \Z$ and $s_\alpha = 1$
if $\alpha \in \Z$, and $Sp(g) = \sum_{(\alpha,w) \in Spp(g)} (\alpha)$ be the spectrum of $g$. 

\begin{rem}
Note that $\alpha$ is a spectral number, i.e. $\exp(-2\pi i\alpha)$ is an eigenvalue of the monodromy. 
$(V(g),0)$ is reducible if and only if $0$ is a spectral number (cf. \cite{DS}). The spectral numbers are 
situated in the interval $(-1,1)$ and the spectrum is symmetric ($\alpha$ is a spectral number 
if and only if $-\alpha$ is a spectral number, cf. \cite{K}). If the Newton polygon of $g$ is non--degenerate
the spectral pairs can be computed combinatorially using the Newton polygon (cf. \cite{Sa}, \cite{St2}).
There is a formula to compute the spectral pairs from the data of the resolution (cf. \cite{SSS}).
\end{rem}

\begin{prop} \label{propSpectralPairs}
$\sigma(z^N+g) = \eta(g,N) - N \cdot \eta(g,1)$ and 
\begin{align*}
\eta(g,K) = \sum_{\substack{\alpha \neq 0, K\alpha \in \Z\\(\alpha,2) \in Spp(g)}} h^{11}_{\exp(-2\pi i\alpha)}
- 2 \sum_{\substack{\alpha \geq 0, K\alpha \notin \Z \\ (\alpha,w) \in Spp(g)}} 
h^{1+[-\alpha], w+s_\alpha-1-[-\alpha]}_{\exp(-2\pi i \alpha)} (1-2 \{K\alpha\}).
\end{align*}
\end{prop}

\begin{rem}
The \textsc{Singular} implementation of the third approach bases on the procedure \texttt{sppairs} (cf.
library \texttt{gmssing.lib}) to obtain the spectral pairs of $(V(g),0)$.  The command \texttt{list L = sppairs(g);} 
provides a list \texttt{L} of the following format: \texttt{L[1]}: set of spectral numbers $\{\alpha_1, \ldots, \alpha_r\}$, 
\texttt{L[2]}: set of weights $\{w_1,\ldots,w_r\}$, \texttt{L[3]}: set of multiplicities $\{h_1,\ldots,h_r\}$ such that 
$Spp(g) = \sum_{j=1}^r h_j \cdot (\alpha_j,w_j)$. The required combinatorics of Proposition 
\ref{propSpectralPairs} had to be implemented in spring 2011.
\end{rem}

\subsection{Theoretical Comparison}

The topological type of a plane curve singularity defined by $g(x,y)=0$ can be described by the
\emph{Puiseux pairs} (Approach 1) of the branches and their \emph{intersection multiplicities} or equivalently 
by discrete \emph{invariants of the resolution} (Approach 2) There are combinatorial formulas to get from one 
description to another (cf. \cite{JP}).
Moreover, the \emph{spectral pairs} (Approach 3) which are topological invariants introduced by Arnold 
(cf. \cite{AGZV}) and Steenbrink (cf. \cite{St2}) can also be computed combinatorially from the resolution data 
(cf. \cite{SSS}). 
Consequently, all of the three approaches to compute the signature of the surface $z^N+g(x,y)=0$ as 
described above are based on the knowledge of three different finite sets of invariants which are related in 
a combinatorial way. The essential difference concerning  these approaches is the method to compute the 
set of the corresponding invariants. 

The \emph{spectral pairs} (Approach 3) can be computed from the mixed Hodge structure. This requires 
several standard basis computations of certain modules over local rings which is the bottleneck of this 
approach. The \textsc{Singular} library \linebreak \texttt{gmssing.lib} is designed for computing the mixed Hodge 
structure for hypersurface singularities of any dimension. This is one reason why computing the 
\emph{spectral pairs} using this method is usually comparatively slow.
 
The \emph{Puiseux pairs} (Approach 1) of the branches and their \emph{intersection multiplicities} resp. 
the \emph{resolution graph} (Approach 2) and the \emph{multiplicity sequence} can be computed via 
Hamburger--Noether expansion (cf. \cite{C}) or resolution of the curve singularities. Hence, both 
approaches are similarly time--consuming. They only need Gr\"obner basis computations if field 
extensions of $\Q$ are necessary to compute the Hamburger--Noether expansion resp. the resolution.  
Anyway, the field extensions are the bottlenecks of these approaches.

\pagebreak

\section{Examples and timings} \label{secExTime}

In this section we provide examples on which we time the three approaches as described in section 
\ref{secTheory} to compute the signature of a surface singularity $z^N+g(x,y) = 0$. 
The corresponding procedures are implemented in \singular in the library \texttt{surfacesignature.lib}.
Timings are conducted by using \singular{3-1-3} on an Intel\textregistered \ Xeon\textregistered \ X5460
with $4$ CPUs,  3.16 GHz each, 64 GB RAM under the Gentoo Linux operating system. 

\begin{exmp} \label{exPoly}
We consider the following polynomials:
\begin{align*}
g_1 = & \;x^{15}-21x^{14}+8x^{13}y-6x^{13}-16x^{12}y+20x^{11}y^2-x^{12}+8x^{11}y-36x^{10}y^2\\
            & \;+24x^9y^3+4x^9y^2-16x^8y^3+26x^7y^4-6x^6y^4+8x^5y^5+4x^3y^6-y^8, \\
g_2 = & \;g_1^3 + x^{17}y^{17}, \\
g_3 = & \;\left(y^4+2x^3y^2+x^6+x^5y\right)^3 + x^{17}y^{17}, \\
g_4 = & \;g_1^5 + x^{20}y^{20}, \\
g_5 = & \;x^{10} + 7y^{10}, \\
g_6 = & \;x^{20} + 5y^{20}.
\end{align*}

The curve singularities in $(\C^2,0)$ defined by $g_1$ resp. $g_2$ resp. $g_3$ are analytically irreducible
with Puiseux pairs $(3,2),(7,2),(15,2)$ resp. $(3,2),(7,2),(15,2),(67,3)$ resp. $(3,2),(7,2),(113,3)$.

The curve singularities in $(\C^2,0)$ defined by $g_4$ resp. $g_5$ resp. $g_6$ are analytically reducible
since they are intersections of $40$ resp. $10$ resp. $20$ lines at the origin. Consequently, the first approach 
is not applicable for these examples. Furthermore, the polynomials $g_4,g_5,g_6$ are defined over $\Q$, 
whereas the resolution is only defined in field extensions of degree $8$ resp. $10$ resp. $20$ over $\Q$. 
\end{exmp} 
 
Computations reveal the following results and corresponding timings which are summarized in Table 
\ref{tabEx1}. The symbol "$> 14\,$h" indicates that the computation did not terminate after more than $14$ 
hours. All timings are given in seconds (s).

\begin{table}[hbt]
\begin{center}
\begin{tabular}{|r|r|r|r|r|r|}
\hline
$N$ & $g_i$ & $\sigma(z^N+g_i(x,y))$ & Approach $1$ [s] & Approach $2$ [s] & Approach $3$ [s] \\
\hline \hline
5 & $g_1$ & $-168$ & 0 & 0 &  $> 14\,$h \\ \hline
5 & $g_2$ & $-1620$ & 174 &  183 & $> 14\,$h \\ \hline
6 & $g_3$ & $-940$ & 2908 & 2912 & $> 14\,$h \\ \hline
5 & $g_4$ & $-3192$ & -- & 19 & $> 14\,$h \\ \hline 
6 & $g_5$ & $-189$ & -- & 22 & 0 \\ \hline
6 & $g_6$ & $-779$ & -- & 14542 & 8 \\ \hline
\end{tabular}
\end{center}
\hspace{15mm}
\caption{Results and total running times for computing the signature of the surface singularity given by
the considered examples (cf. Example \ref{exPoly}) via all approaches as described in section \ref{secTheory}. } 
\label{tabEx1}
\end{table}

In addition, we summarize the maximal memory allocated from system during the considered computations 
in Table \ref{tabEx2}. The memory consumption is given in Megabyte (MB).  

\begin{table}[hbt]
\begin{center}
\begin{tabular}{|r|r|r|r|r|r|}
\hline
$N$ & $g_i$ & Approach $1$ [MB] & Approach $2$ [MB] & Approach $3$ [MB] \\
\hline \hline
5 & $g_1$ & $1$ & $1$ &  $> 300$ \\ \hline
5 & $g_2$ & $3442$ & $3442$ & $> 1100$ \\ \hline
6 & $g_3$ & $7723$ & $7723$ & $> 1100$ \\ \hline
5 & $g_4$ & -- & $33$ & $> 2300$ \\ \hline 
6 & $g_5$ & -- & $8$ & $4$ \\ \hline
6 & $g_6$ & -- & $156$ & $72$ \\ \hline
\end{tabular}
\end{center}
\hspace{15mm}
\caption{Maximal memory allocated from system while computing the signature of the surface singularity 
given by the considered examples (cf. Example \ref{exPoly}) via all approaches as described in section 
\ref{secTheory}. } 
\label{tabEx2}
\end{table}

\begin{rem}[Algorithmic Conclusion]
Our experiments reveal that there exist examples where Approach $1$ and Approach $2$ are almost
equivalent regarding time consumption resp. memory allocation, but Approach $3$ is much slower
(cf. $g_1$, $g_2$, $g_3$). Furthermore, there exist examples where Approach $1$ is not applicable, but
Approach $2$ consumes more time resp. allocates more memory than Approach $3$ (cf. $g_4$), and
vice versa (cf. $g_5$, $g_6$).
Consequently, it is reasonable to summarize all approaches in one algorithm which computes the 
signature via every approach, if possible, in parallel such that the fastest approach wins and returns the 
result. 
\end{rem}

\end{document}